\documentclass{icmart}

\contact[oguiso@math.sci.osaka-u.ac.jp]{Department of Mathematics, Osaka University, Toyonaka 560-0043, Osaka, Japan}

\theoremstyle{plain}
\newtheorem{thm}{Theorem}[section]
\newtheorem{theorem}[thm]{Theorem}

\newtheorem{corollary}[thm]{Corollary}

\theoremstyle{definition}
\newtheorem{remark}[thm]{Remark}

\newtheorem{definition}[thm]{Definition}

\newtheorem{example}[thm]{Example}

\newtheorem{question}[thm]{Question}

\newtheorem{problem}[thm]{Problem}

\numberwithin{equation}{section}

\title [Explicit birational geometry inspired by complex dynamics]
{Some aspects of explicit birational geometry inspired by complex dynamics}
\author[Keiji Oguiso]
{Keiji Oguiso}

\begin{document}

\maketitle

\begin{abstract} 
Our aim is to illustrate how one can effectively apply the basic ideas and notions of topological entropy and dynamical degrees, together with recent progress of minimal model theory in higher dimension, for an explicit study of birational or biregular selfmaps of projective or compact K\"ahler manifolds, through 
concrete examples. 
\end{abstract}

\section{Introduction}

This is a survey of some aspects of recent progress on birational and biregular complex algebraic geometry inspired by complex dynamics in several variables. 
Our aim is to illustrate how one can apply the basic ideas and notions of topological entropy and dynamical degrees, together with recent progress of minimal model theory in higher dimension, for an explicit study of birational or biregular selfmaps of projective or compact K\"ahler manifolds. Especially, we focus on the following one of the most basic, natural problems:

\begin{problem}\label{Explicit1}
Find many examples of projective or K\"ahler manifolds $M$ admitting {\it interesting} birational automorphisms of {\it infinite order}, 
or more preferably, {\it primitive biregular} automorphisms 
of {\it positive entropy}. 
\end{problem} 

There are so many interesting works in this area since the breakthrough results due to Serge Cantat \cite{Ca99} and Curtis T. McMullen \cite{Mc02-2}, and this 
note is definitely far from being a complete panorama of this area. Also, needless to say, there is no universally acceptable mathematical definition of the term {\it interesting} and the choice of topics and materials owes much to my own flavour and ability, and probably not the one that everyone agrees with. For instance, the terms {\it of infinite order} ignore very beautiful aspects of finite group actions on manifolds. 

Throughout this note, we work over the complex number field ${\mathbf C}$. We 
assume some familiarity with basics on complex geometry and algebraic geometry. Unless stated otherwise, the topology we use is the Euclidean topology (not Zariski topology), a point means a closed point, and manifolds and varieties are connected. By abuse of language, we call a (bi)meromorphic map also a (bi)rational map even under non-algebraic settings. We denote by $I(f)$ the indeterminacy locus of a rational map $f : M \cdots\to N$, i.e., the complement of the maximal, necessarily Zariski open dense, subset $U \subset M$ such that $f \vert U$ is holomorphic. $I(f)$ is a Zariski closed subset of $M$ of codimension $\ge 2$ if $M$ is normal. The set of birational selfmaps (resp. biregular selfmaps, resp. birational selfmaps being isomorphic in codimension one) of $M$ form a group under the natural composition. We denote these groups by ${\rm Bir}\, (M)$ (resp. ${\rm Aut}\, (M)$, resp. ${\rm PsAut}\, (M)$). Here we call $f$ isomorphic in codimension one if $f$ neither contracts nor extracts any divisors, in other words, if there are Zariski open dense subsets $U$, $U'$ of $M$ such that ${\rm codim}\, M \setminus U \ge 2$, ${\rm codim}\, M \setminus U' \ge 2$ for which the restriction map $f\vert U : U \to U'$ is an isomorphism. We call an element of ${\rm Bir}\, (M)$ (resp. ${\rm PsAut}\, (M)$, resp. ${\rm Aut}\, (M)$) a birational automorphism (resp. a pseudo-automorphism, resp. an automorphism or a biregular automorphism) of $M$. We call a compact complex variety of class ${\mathcal C}$ if it is birational to a compact K\"ahler manifold. Unless stated otherwise, 
we denote 
the golden number $(\sqrt{5} +1)/2$ by $\eta$, the cyclic group of order $n$ by ${\mathbf Z}_n$, the free product of groups $G_1$ and $G_2$ by $G_1 * G_2$.

\section{What kind of manifolds we are interested in?}

Let $M$ be a compact K\"ahler manifold of $\dim M = l >0$. We are interested 
in a {\it birational automorphism} $f$ of $M$, in particular, 
{\it of infinite order}.

\subsection{Primitive birational automorphisms after De-Qi Zhang} If $f_i \in {\rm Bir}\, (M_i)$ ($i =1$, $2$), then $f_1 \times f_2 \in {\rm Bir}\, (M_1 \times M_2)$ and it is of infinite order if so is one of $f_i$. We are more interested in birational automorphisms {\it not coming from lower dimensional pieces}, more precisely, {\it primitive} ones in the sense of 
De-Qi Zhang \cite{Zh09-1}:

\begin{definition}\label{primitive} A birational automorphism $f$ of $M$ is {\it imprimitive} if there are a dominant rational map $\varphi : M \cdots \to B$ to a compact complex variety $B$ with $0 <\dim\, B < \dim\, M = l$ and a rational map $g : B \cdots \to B$, necessarily a birational automorphism of $B$, such that $\varphi \circ f = g \circ \varphi$. A birational automorphism that is not imprimitive is {\it primitive}. 
\end{definition}

Here we may assume that $B$ is smooth and $\varphi : M \to B$ is holomorphic whenever it is more convenient. Indeed, we may resolve $B$ first and then resolve the indeternimacy of $\varphi$, by the fundamental result of Hironaka. 

{\it What kind of manifolds can have primitive birational automorphisms of infinite order?} We can not answer this question completely, but if we {\it assume} that minimal model program (\cite{KMM87}, \cite{KM98}, \cite{HM10}) in higher dimensional projective manifolds work, then one has the following rough but quite nice picture at least in the projective case. This beautiful observation is due to De-Qi Zhang (see \cite{Zh09-1} and \cite{NZ09} for more precise results):

\begin{theorem}\label{deqizhang} Let $M$ be a projective manifold of dimension $l$. Assume that the minimal model conjecture (MMP) and the weak abundance (WA) in dimension $l$ hold, in the sense that any $l$-dimensional projective manifold $V$ is birational to either a minimal model with weak abundance or a Mori fiber space, i.e., there is a projective variety 
$V_{{\rm min}}$, birational to $V$,  with only normal, ${\mathbf Q}$-factorial terminal singularities 
such that either one of the following two holds:

(MMP-1) in addition, $K_{V_{{\rm min}}}$ is nef (minimal model); or 

(MMP-2) in addition, there exists an extremal contraction $\varphi : V_{{\rm min}} \to B$, with respect to $K_{V_{{\rm min}}}$ onto a normal projective variety such that $0 \le \dim\, B < \dim\, V$ (Mori fiber space), 

and additionally,

(WA) In the case (MMP-1), $\vert mK_{V_{{\rm min}}} \vert$ 
is non-empty for some $m>0$.

Then any $l$-dimensional projective manifold $M$ with a primitive birational automorphism $f$ of infinite order is birational to either:

(RC) a rationally connected manifold, in the sense that any two points can be connected by a finite chain of rational curves;

(WCY) a minimal Calabi-Yau variety, in the sense that it is a minimal variety with numerically trivial canonical divisor and of irregularity $0$; or

(T) an abelian variety, i.e., a projective complex torus. 

\end{theorem}

\begin{example}\label{lowdimdeqizhang}
(1) When $l = 1$, (RC) is ${\mathbf P}^1$, (T) is an elliptic curve and no (WCY). 

(2) When $l=2$, (RC) is a rational surface, (T) is an abelian surface, and (WCY) is a K3 surfaces or an Enriques surface (see eg. \cite{BHPV04}).
\end{example}

Since the proof of Theorem (\ref{deqizhang}) provides a good introduction on objects we are interested in, we give here a fairly complete proof, following \cite{Zh09-1}, \cite{Ue75}. 

\begin{proof} The most essential point is that any {\it canonically defined maps} are preserved by ${\rm Bir}\, (M)$. 

Let $\kappa\, (M) \in \{-\infty, 0, 1, \cdots, l\}$ be the Kodaira dimension. $\kappa\, (M)$ is the maximal dimension of the images $W_m$ under the pluri-canonical maps associated to the complete linear system $\vert mK_M \vert$ ($m = 1, 2, 3, \cdots $) :
$$\Phi_m := \Phi_{\vert mK_M \vert} : M \cdots\to W_m := {\rm Im}\,\Phi_{\vert mK_M \vert} \subset \vert mK_M \vert^* = {\mathbf P}^{\dim\, \vert mK_M \vert}\,\, ,$$
if $\vert mK_M \vert \not= \emptyset$ for some $m > 0$ and $\kappa\, (M) = -\infty$ otherwise. It is a birational invariant. 

Consider first the case where $\kappa (M) \ge 0$. We may assume that $\Phi_m$ is regular. We can and do choose $m$ so that $\dim\, W_m = \kappa (M)$ and $\Phi_m$ is of connected fibers. We write $W = W_{m}$. Note that any birational automorphism 
$f$ preserves the set of global holomorphic pluri-canonical forms, as ${\rm codim}\, I(f) \ge 2$. Then the induced projective linear map $f_* \in {\rm PGL}\, (\vert mK_M \vert^{*}) = {\rm Aut}\, (\vert mK_M \vert^*)$ preserves $W$ and is equivariant to $f$ with respect to $\Phi_m$. Hence $f$ is imprimitive if $1 \le \kappa (M) \le l -1$. 

If $\kappa (M) = l$, then the same is true but $M$ and $W$ are birational. So ${\rm Bir}\, (M) = {\rm Aut}\, (W) \subset {\rm PGL}\, (\vert mK_M \vert^{*})$. It is Zariski closed in the affine noetherian group ${\rm PGL}\, (\vert mK_M \vert^*)$, as it is the stabilizer of the point $[W]$ of the action of ${\rm PGL}\, (\vert mK_M \vert^{*})$ on ${\rm Hilb}\, (\vert mK_M \vert^*)$. Hence it is finite. Indeed, if otherwise, $\dim\, {\rm Aut}\, (W) \ge 1$ and we can choose a one dimensional algebraic subgroup, which is necessarily isomorphic to ${\mathbf C}$ or ${\mathbf C}^{\times}$. The Zariski closures of the orbits of general points of $W$ under this $1$-dimensional algebraic subgroup are necessarily rational, and cover $W$ and $M$. Then again by the Hilbert scheme, we have a dominant holomorphic maps $\pi : Y \to M$ from the fiber space $Y \to X$ 
whose general fibers $Y_x$ are isomorphic to ${\mathbf P}^1$. Since we work over the field of characteristic $0$, the map $\pi$ is separable so that $\vert mK_Y \vert \not= \emptyset$ as well. Hence for general $Y_x$, we have ${\rm deg}\, mK_{Y_x} \ge 0$ by the adjunction formula, a contradiction to the fact that $Y_x \simeq {\mathbf P}^1$. Hence ${\rm Bir}\, (M)$ is finite and in particular, $f$ is of finite order when $\kappa (M) = l$. 

Hence $\kappa (M) = 0$ or $-\infty$. So far we did not use the assumption (MMP), (WA). Also, our assumption that $f$ is of infinite order is used only to conclude $\kappa (M) \not= l$. 

Assume that $\kappa (M) = 0$. Consider next the irregularity $q(M) := h^0(M, \Omega_M^1)$.  If $q(M) > 0$, then we have the albanese morphism
$${\rm alb}_M : M \to {\rm Alb}\, (M) = H^0(M, \Omega_M^1)^{*}/H_1(M, {\mathbf Z})\,\, .$$
It is classical that ${\rm Alb}\, (M)$ is an abelian variety. Since $\kappa (M) = 0$, a fundamental theorem due to Kawamata \cite{Ka81} (again free from (MMP), (WA)) says that ${\rm alb}_M$ is surjective with connected fibers, in particular $q(M) \le l$. For the same reason as before, the action ${\rm Bir}\, (M)$ descends to the biregular action of ${\rm Alb}\, (M)$ equivariantly with respect to the albanese map. Hence either $q(M) = 0$ or $q(M) = l$. In the second case, $M$ is birational to ${\rm Alb}\, (M)$. If $q(M) = 0$, then by our assumption (MMP), $M$ is birational to a minimal Calabi-Yau variety. (Here, if one prefers, one also stops at the stage $\kappa (M) = 0$ and $q(M) = 0$. Then the conjectural (MMP) is not required here.) 

It remains to treat the case $\kappa (M) = -\infty$. This is the most subtle case where we really use our assumptions (MMP) and (WA) (but we do not use our assumption that $f$ is of infinite order anymore). (WA) is one way to conclude (MMP-1) and (MMP-2) are exclusive. Let us consider $W := M_{{\rm min}}$ in Theorem (\ref{deqizhang}), whose existence needs (MMP). If $K_W$ would be nef, then by (WA), $\vert mK_M \vert = \vert mK_W \vert \not= \emptyset$ for some $m >0$, a contradiction to $\kappa (M) = -\infty$. Hence, the case (MMP-1) does not happen and therefore (MMP-2) happens by our assumptions (MMP). In (MMP-2), by the property of an extremal contraction, the fibers of the Mori fiber space are covered by rational curves. Now we consider the maximal rationally connected fibration, MRC fibration, for short (\cite{Ko96}). The MRC fibration $r : M \cdots\to R$ is an almost holomorphic, rational dominant map such that for general $x \in M$, the fiber $M_p \ni x$ is the maximal rationally connected submanifold of $M$ containing $x$, and birationally preserved by ${\rm Bir}\, (M)$, hence by $f$. Since $0 \le \dim\, R \le l -1$ for our $M$, it follows that $R$ is a point, i.e., 
$M$ is rationally connected. 
\end{proof} 
(MMP) and (WA) hold in dimension $\le 3$, finally due to Mori \cite{Mo88} (MMP) and Kawamata \cite{Ka92} (WA) in the strongest form that nef $K_{V_{{\rm min}}}$ is actually semi-ample. So, Theorem (\ref{deqizhang}) is unconditional in dimension $\le 3$. In higher dimension, both conjectures are expected to be true (cf. \cite{BCHM10}, \cite{HM10}, \cite{Na04}). In dimension $2$, Theorem (\ref{deqizhang}) is essentially the same as the breakthough observation due to Cantat \cite{Ca99} (see also \cite{CCG10}), in terms of {\it topological entropy}. 

\subsection{Three classes of manifolds in Theorem (\ref{deqizhang})} In this subsection, we discuss basic manifolds belonging to the three classes in Theorem (\ref{deqizhang}), which are indeed the main objects in this note.

{\bf Rational manifolds.} An excellent reference of rationally connected manifolds (RC manifolds) is \cite{Ko96}. Most basic examples of RC manifolds are rational manifolds, i.e., manifolds which are birational to ${\mathbf P}^l$. Note that ${\rm Aut}\, ({\mathbf P}^l) = {\rm PGL}(l+1, {\mathbf C})$. It is obvious that generic $g \in {\rm Aut}\, ({\mathbf P}^l)$ is of infinite order. On the other hand, any $g \in {\rm Aut}\, ({\mathbf P}^l)$ is imprimitive if $l \ge 2$. Indeed, $g$ has a fixed point $P \in {\mathbf P}^l$, corresponding to the eigenvector of a lift $\tilde{g} \in {\rm GL}(l+1, {\mathbf C})$. The family of lines through $P$ is then stable under $g$. Let ${\rm Bl}_P\, {\mathbf P}^{l}$ be the blow up of ${\mathbf P}^l$ at $P$. The lines through $P$ are the fibers of the natural morphism ${\rm Bl}_P\, {\mathbf P}^{l} \to {\mathbf P}^{l-1} = {\mathbf P}(T_{P, {\mathbf P}^l})$, and this fibration is stable under the natural (biregular) action of $g$. 

So, {\it in our view}, ${\rm Aut}\, ({\mathbf P}^l) = {\rm PGL}\, (l+1, {\mathbf C})$ is not so interesting. However, the group ${\rm PGL}\, (l+1, {\mathbf C})$ has a very deep aspect in birational geometry, for instance, the following striking result due to Cantat \cite{Ca13}: 
\begin{theorem}\label{cantat1} Let $M$ be an $l$-dimensional projective 
manifold. If there is an injective group homomorphism ${\rm PGL}\, (n+1, {\mathbf C}) \to {\rm Bir}\, (M)$, as abstract groups, then $n \le l$ and the equality $n = l$ holds if and only if $M$ is rational. In particular, ${\rm Bir}\, ({\mathbf P}^{l}) \simeq {\rm Bir}\, ({\mathbf P}^{l'})$ as abstract groups if and only if $l = l'$. 
\end{theorem}

The {\it standard Cremona transformation} 
$${\rm cr}_l : {\mathbf P}^l \cdots\to {\mathbf P}^l\,\, ,\,\, [x_0 : x_1 : \cdots : x_{l}] \mapsto [\frac{1}{x_0} : \frac{1}{x_1} : \cdots : \frac{1}{x_l}]$$
is the most basic birational non-biregular automorphism of ${\mathbf P}^l$ ($l \ge 2$). 
The indeterminacy locus of ${\rm cr}_l$ are $\cup_{i \not= j} L_{ij}$, where $L_{ij} := (x_i = x_j = 0)$. Let ${\rm SCr}_{l} := \langle {\rm PGL}(l+1, {\mathbf C}), {\rm cr}_l \rangle < {\rm Cr}_l := {\rm Bir}\, ({\mathbf P}^{l})$. 
We have ${\rm SCr}_2 = {\rm Cr}_2$ (Noether's theorem, \cite{Do11}). If $l \ge 3$, then ${\rm SCr}_l$ is 
much smaller than ${\rm Cr}_l$ but ${\rm SCr}_l$ is rich enough. One of unexpected applications of ${\rm SCr}_l$ 
is the following result due to Lesieutre 
\cite{Le13}. In \cite{Le13}, the group ${\rm SCr}_3$ and its complex dynamical aspect are effectively applied to prove the following derived categorical result:

\begin{theorem}\label{lesieutre} There is a smooth rational threefold with infinitely many birational non-isomorphic Fourier-Mukai partners. 
\end{theorem}

Note that ${\rm cr}_l$ maps the coordinate hyperplane $H_i = (x_i = 0)$ to the standard coordinate point $e_i = [0 : \cdots : 0 : 1 : 0 : \cdots 0]$, where $1$ is at the $i$th coordinate. So, ${\rm cr}_l \in {\rm Bir}\, ({\mathbf P}^l) \setminus {\rm PsAut}\, ({\mathbf P}^l)$. Actually ${\rm PsAut}\, ({\mathbf P}^l) = {\rm Aut}\, ({\mathbf P}^l)$ by ${\rm Pic}\, {\mathbf P}^l = {\mathbf Z}H$. However, if we blow-up ${\mathbf P}^l$ at the $(l+1)$ standard coordinate points $e_i$, then ${\rm cr}_l$ gives rise to a pseudo-automorphism 
$\tilde{{\rm cr}}_i$ of
${\rm Bl}_{\{e_i\}} {\mathbf P}^l$, and performing further blowing-ups, we can make it a biregular automorphism, of order $2$. Let $S, T \in {\rm PGL}(l+1, {\mathbf C})$. Then $f = S \circ {\rm cr}_l \circ T^{-1} \in 
{\rm SCr}_l$ is of infinite order for almost all choices of 
$S$ and $T$, and $f$ lifts to a {\it pseudo-automorphism} of some blowing-ups of ${\mathbf P}^l$, under some periodicity condition for the indeterminacy loci $I(f^{\pm n})$ (\cite{BK06}, \cite{BK09}, \cite{BK14}). In this way, Bedford and Kim construct many interesting rational surface automorphisms as well as pseudo-automorphisms of rational threefolds very explicitly. However, when $l =3$, none of them seems to be realized as a {\it biregular} automorphism (cf. Question (\ref{bedford})).  

{\bf A few properties of pseudo-automorphisms.} Before entering two other classes, we recall a few basic properties of ${\rm PsAut}\, (M)$. The group ${\rm PsAut}\, (M)$ naturally acts on the N\'eron-Severi group ${\rm NS}\, (M) := {\rm Im}(c_1 : {\rm Pic}\, (M) \to H^2(M, {\mathbf Z}))$ as well as on $H^2(M, {\mathbf Z})$. This action is functorial, in the sense that $(f \circ g)^* = g^* \circ f^*$ on $H^2(M, {\mathbf Z})$ and preserves the Hodge decomposition of $H^2(M, {\mathbf Z})$ (but not the intersection $(x^l)$ in general). For a minimal model in the sense (MMP-1), we have the following factorization. This fundamental result is due to Kawamata \cite{Ka08}: 
 
\begin{theorem}\label{kawamata1} Let $M$ be a minimal model. Then ${\rm Bir}\, (M) = {\rm PsAut}\, (M)$, and any $f \in {\rm Bir}\, (M)$ is decomposed as $f = \varphi \circ \iota_{m-1} \circ \cdots \circ \iota_0$, where $M_0 = M = M_m$ and $\iota_i : M_i \cdots\to M_{i+1}$ ($0 \le i \le m-1$) are flops betweeen minimal models $M_i$ and $M_{i+1}$ and $\varphi \in {\rm Aut}\, (M)$. 
\end{theorem}

{\bf Complex tori, CY manifolds and HK manifolds.} The case of complex tori is very much known (See \cite{FZ13}, \cite{OT13-1} for some dynamically interesting features in tori). Most basic examples of minimal Calabi-Yau varieties are CY manifolds and (projective) HK manifolds as defined below. CY manifolds, HK manifolds together with rational manifolds are the main objects in this note. 

\begin{definition}\label{cyhk} Let $M$ be an $l$-dimensional simply-connected compact K\"ahler manifold. 

(1) $M$ is a Calabi-Yau manifold in the strict sense ({\it CY manifold}) if $l \ge 3$ and 
$H^0(M, \Omega_M^j) = 0$ for $0 < j \le l-1$ and $H^0(M, \Omega_M^l) = {\mathbf C}\omega_M$, where $\omega_M$ is a nowhere vanishing holomorphic $l$-form.

(2) $M$ is a compact hyperk\"ahler manifold ({\it HK manifold}) if $H^0(M, \Omega_M^2) = {\mathbf C}\sigma_M$, where $\sigma_M$ is an everywhere nondegenerate holomorphic $2$-form.
\end{definition}

Good references of CY manifolds and HK manifolds are \cite{GHJ03}, \cite{Ma11}. By definition, HK manifolds are of even dimension and K3 surfaces are nothing but HK manifolds of dimension $2$. CY manifolds $M$ are always projective by $h^0(\Omega_M^2) = 0$ (as $l \ge 3$). On the other hands, both projective HK manifolds and non-projective HK manifolds are dense both in the Kuranishi space and in the global moduli space of marked HK manifolds (\cite{Fu83}, \cite{Hu03}). Examples with interesting (birational) automorphisms in our view will be given in Sections 4, 5, 6. 

The importance of CY manifolds and HK manifolds in complex algebraic geometry lies in the fact, called the Bogomolov decomposition theorem (\cite{Be84}), that 
these two classes of manifolds together with complex tori form the building blocks of compact K\"ahler manifolds with trivial first Chern class.

We close this section by the following:

\begin{remark}\label{ample} Let $M$ be a CY manifold or a projective 
HK manifold and $G < {\rm Bir}\, (M) = {\rm PsAut}\, (M)$. Assume that 
there is an ample divisor $H$ such that  
$f^*H = H$ in ${\rm Pic}\, (M) \simeq {\rm NS}\, (M)$ for all $f \in G$. 
Then $G < {\rm Aut}\, (M)$ 
and $G$ is a finite group. In particular, if $\rho (M) := {\rm rank}\, {\rm NS}\, (M) =1$, then ${\rm Bir}\, (M) = {\rm Aut}\, (M)$ and it is a finite group. {\it So, in our view, interesting cases are $\rho (M) \ge 2$.} 
\end{remark}
Indeed, the same argument as in Theorem (\ref{deqizhang}), applied for the 
$G$-equivariant embedding $\Phi_{\vert mH \vert} : M \to \vert mH \vert^*$ for large $m >0$, shows that 
$G < {\rm Aut}\, (M)$ and at the same time $G$ is a Zariski closed 
algebraic subgroup of ${\rm PGL}\, (\vert mH \vert^*)$. Since $\dim\, G = 0$ by $H^0(M, TM) = 0$, the result follows.

\section{Topological entropy and Dynamical degrees}

\subsection{Topological entropy} References of this subsection are 
\cite{KH95}, \cite{Gr87}.

Let $X = (X, d)$ be a compact metric space and $f : X \rightarrow X$ be a continuous surjective selfmap of $X$. We denote by $f^n$ the $n$-th iterate of $f$. The {\it topological entropy} of $f$ is the fundamental invariant that measures {\it how fast two general points spread out under the action of the semi-group $\{f^n \vert n \in {\mathbf Z}_{\ge 0}\}$}, hence, presents a kind of complexity of $f$. For the definition, we define the new distance $d_{f, n}$ on $X$ by 
$$d_{f, n}(x, y) = {\rm max}_{0 \le j \le n-1} d(f^j(x), f^j(y))\,\, {\rm for}\,\, x, y \in X.$$
Under the identification $x \leftrightarrow {\mathbf x}^{(n)} := (x, f(x), \cdots , f^{n-1}(x))$, 
the new distance $d_{f, n}(x, y)$ is the distance of the graph
$$\Gamma_{f, n} := \{{\mathbf x}^{(n)} = (x, f(x), \cdots , f^{n-1}(x))\, 
\vert\, x \in X\} \subset X^n$$
induced by the product distance on $X^d$. The first projection ${\rm pr}_1 : (\Gamma_{f, n}, d_{f, n}) \to (X, d_{f, n})$ is an isometry and ${\rm pr}_1 : (\Gamma_{f, n}, d_{f, n}) \to (X, d)$ is a homeomorphism.  
 
Let $\epsilon > 0$ be a positive real number. 
We call two points $x, y \in X$ $(n, \epsilon)$-{\it separated} if $d_{f, n}(y, x) \ge \epsilon$, and a subset $F \subset X$ $(n, \epsilon)$-{\it separated} if any two distinct points of $F$ are $(n, \epsilon)$-separated. 
Let 
$$N_{d}(f, n, \epsilon) := {\rm Max}\, \{\vert F \vert\,\, \vert\,\, F \subset X\,\, {\rm is}\,\,  (n, \epsilon)-{\rm separated}\,\, \}\,\, .$$
Note that $N_{d}(f, n, \epsilon)$ is a well-defined positive integer, because $X$ is compact. 

\begin{remark}\label{eye} Please imagine that $\epsilon >0$ is "very very 
small", so that we can {\it not} distinguish two points $x, y \in X$ with $d(x, y) < \epsilon$ by "our eyes" but can do if $d(x, y) \ge \epsilon$. Then, we can not distinguish $x, y$ if they are not $(1, \epsilon)$-separated but we can distinguish them by performing $f$ if they are $(2, \epsilon)$-separated. Similarly, we can distinguish $x, y$ at some stage, say $f^j(x), f^j(y)$ ($0 \le j \le n-1$), if they are $(n, \epsilon)$-separated. In this sense, $N_{d}(f, 1, \epsilon)$ is the maximal number of points of $X$ distinguished by eyes and $N_{d}(f, n, \epsilon)$ is the maximal number of points of $X$ distinguished by eyes after performing $f^j$ ($0 \le j \le n-1)$. So, roughly, the growth of the sequence $\{N_{d}(f, n, \epsilon)\}_{n \ge 1}$ measures how fast general points spread out under the iterations of $f$ to be distinguishable by our eyes (if it will be). 
\end{remark} 

\begin{definition}\label{bowen} The {\it topological entropy}, or {\it entropy} for short, of $f$ is:
$$h_{{\rm top}}(f) := h_d(f) := {\rm lim}_{\epsilon \to +0} h_d (f, \epsilon)\,\, ,$$
where
$$h_d (f, \epsilon) := {\rm limsup}_{n \to \infty} \frac{\log N_d(f, n, \epsilon)}{n}\,\, .$$
\end{definition}
    
Since $\log N_d(f, n, \epsilon) \ge 0$ is an increasing function of $\epsilon > 0$, the limit exists in $[0, \infty]$ (possibly $\infty$). More or less from the definition, we obtain:

\begin{corollary}\label{entropybasic}

(1) $h_{{\rm top}}(f)$ is a topological invariant, in the sense that $h_{d'} (f) = h_d (f)$ for any distance $d'$ of $X$ such that $(X, d')$ and $(X, d)$ are homeomorphic. 

(2) If $h_{{\rm top}}(f) > 0$, then $f^m \not= id_X$ for all $m \ge 1$, i.e., ${\rm ord}\, (f) = \infty$. 

(3) If $f$ is an isometry, for instance a translation of a torus, then $h_{{\rm top}}(f) = 0$. In particular, the converse of (2) is not necessarily true. 

(4) $h_{{\rm top}}(f \times f') = h_{{\rm top}}(f)+h_{{\rm top}}(f')$ where $f'$ is a surjective selfmap of a compact metric space $X' =(X', d')$. 
\end{corollary} 

\begin{example}\label{tori1} Let $E$ be a $1$-dimensional complex torus. 
Consider the abelian surface $A = E \times E$ and its surjective endomorphism 
$f_M(x) = Mx$ given by $M \in M(2, {\mathbf Z})$ with ${\rm det}\, M \not= 0$. Note that $f_M \in {\rm Aut}\, (A)$ if ${\rm det}\, M = \pm 1$. Let $\alpha$, $\beta$ be the eigenvalues of $M$ such that $\vert \alpha \vert \le \vert \beta \vert$. Then, according to the three cases (i) $\vert \alpha \vert \ge \vert \beta \vert \ge 1$, (ii) $\vert \alpha \vert \ge 1 \ge \vert \beta \vert$, (iii) $1 \ge \vert \alpha \vert \ge \vert \beta \vert$, the entropy $h_{{\rm top}}(f_M)$ 
is (i) $\log \vert \alpha \beta \vert^2$, (ii) $\log \vert \alpha \vert^2$, 
(iii) $\log 1 = 0$. In particular, $h_{{\rm top}}(f_M) = \log \eta^2 >0$, the natural logarithm of (the square of) the golden number, for the Lie automorphism $f_M \in {\rm Aut}\, (A)$ given by 
$$M = \left(\begin{array}{rr} 
2 & 1\\
-1 & -1\\
\end{array} \right)\,\, .$$
\end{example}

Very rough idea is as follows. For simplicity, we further assume that $M$ is diagonalizable in $M(2, {\mathbf C})$. We fix the flat distance $d$ on $A$ from the universal cover ${\mathbf C}^2$. Let $\epsilon > 0$ be a very small number. Let us cover $A$ by $N$ mutually disjoint complex 2-dimensional $\epsilon$-"parallelograms" (actual real dimension is $4$) that are parallel to the complex eigenvectors of $M$. Then $N(f_M, 1, \epsilon)$ is about $N$. Next divide each of $N$ $\epsilon$-parallelograms into mutually disjoint $\epsilon$-parallelograms with respect to the distance $d_{f_M, 2}$. In case (i), each original parallelogram is devided into about $\vert\alpha \beta \vert^2$ new parallelograms, because $\vert \alpha \vert \ge 1$ and $\vert \beta \vert \ge 1$ (and real dimension is $2 +2$). Therefore, $N(f_M, 2, \epsilon)$ is about $\vert\alpha \beta \vert^2 N$. In case (ii), each parallelogram is devided into $\vert\alpha\vert^2$ new parallelogram, because $\vert \alpha \vert \ge 1$ but $\vert \beta \vert \le 1$. Therefore, $N(f_M, 2, \epsilon)$ is about $\vert\alpha\vert^2 N$. Similarly, in case (iii), $N(f_M, 2, \epsilon)$ remains $N$. Repeating this, we see that $N(f_M, n, \epsilon)$ is about $\vert \alpha \beta 
\vert^{2(n-1)}N$, $\vert \alpha \vert^{2(n-1)}N$, $N$ according to the three cases (i), (ii), (iii). This implies the result. 

Note that in each case, the entropy is the natural logarithm of the spectral radius of $f_M^* \vert H^{*}(A, {\mathbf Z})$. Actually, this is {\it not accidental} as we will explain in the next subsection.

\subsection{Fundamental theorem of Gromov-Yomdin} References of this subsection are \cite{Gr87}, \cite{Yo87}, \cite{Gr03}, \cite{DS05-1}, \cite{DS05-2} (see also \cite{CCG10}). 

Let $M$ be a compact K\"ahler manifold of dimension $l$ and $\eta$ be any K\"ahler form on $M$. Then $M$ is a compact metric space by the distance defined by $\eta$. Let $f : M \to M$ be a surjective holomorphic map. Then 
$f^{*}$ naturally acts on the $k$-th cohomology group $H^{k}(M, {\mathbf Z})$ as well as each Hodge component $H^{p,q}(M)$. We define $r_{p}(f)$ to be the {\it spectral radius} of $f^{*} \vert H^{p, p}(M)$, that is, the maximum absolute value of eigenvalues of $f^{*} \vert H^{p, p}(M)$. Similarly, we denote by $r(f)$ (resp. $r^{{\rm even}}(f)$) the spectral radius of $f^{*}$ on $\oplus_{k=0}^{2l} H^{k}(M, {\mathbf Z})$ (resp. $\oplus_{p = 0}^{l} H^{2p}(M, {\mathbf Z})$). 

We define the {\it $p$-th dynamical degree} $d_{p}(f)$ by 
$$d_{p}(f) := \lim_{n \to \infty} (\delta_{p}(f^n))^{\frac{1}{n}}\,\, ,$$ 
where 
$$\delta_{p}(f^n) := (\int_{M} (f^n)^{*}(\eta^{p}) \wedge \eta^{l-p}) = ([(f^n)^{*}(\eta^{p})].[\eta^{l-p}])_M\,\, .$$ 
Here $(*, **)_M$ is the intersection number. The limit does not depend on the choice of $\eta$ once the existence is guaranteed. Indeed, for two K\"ahler forms $\eta$ and $\eta'$, there are positive real number $C$ and a K\"ahler form $\eta"$ such that $C[\eta] = [\eta'] + [\eta"]$ in $H^{1,1}(M, {\mathbf R})$. The fact that the limit exists is non-trivial. There are many ways to see it. For instance, it is an immediate consequence of the following crucial observation by Dinh-Sibony, which holds also for {\it rational} dominant selfmaps (\cite{DS05-1}, \cite{DS05-2}):
\begin{theorem}\label{DinhSibony1}
There is a constant $C = C_{M, \eta}$ depending only on $M$ and 
$\eta$ (but not on $f$ and $g$) such that 
$$\delta_{p}(f\circ g) \le C \delta_{p}(f)\delta_{p}(g)\,\, ,$$ 
for any two dominant holomorphic selfmaps $f : M \to M$, $g : M \to M$. 
\end{theorem}

The logarithmic volume ${\rm lov}\, (f)$, introduced by Gromov, is
$${\rm lov}\, (f) := {\rm limsup}_{n \to \infty} \frac{\log {\rm Volume} (\Gamma_{f, n})}{n}\,\, ,$$
where
$${\rm Volume} (\Gamma_{f, n}) := \frac{1}{l!}\int_{\Gamma_{f, n}} 
(\sum_{i=1}^{n} {\rm pr}_i^*\eta_M)^l\,\, .$$ 
The following fundamental theorem is due to Gromov-Yomdin:
\begin{theorem}\label{GromovYomdin}
Let $M$ be a compact K\"ahler manifold of dimension $l$ and $f : M \to M$ be a surjective holomorphic map. Then, $d_p(f) = r_{p}(f)$ and 
$$h_{{\rm top}}(f) = {\rm lov}\, (f) = \log {\rm max}_{0 \le p \le l}\, d_{p}(f) = \log {\rm max}_{0 \le p \le l}\, r_{p}(f) = \log r^{{\rm even}}(f) = \log r(f)\,\, .$$
Moreover, if $M$ is projective, then $h_{{\rm top}}(f)$ is also equal to the natural logarithm of the spectral radius of $f^* \vert\oplus_{p} H^{p, p}(M, {\mathbf Z})$. 
\end{theorem}

This theorem opens the door to compute the entropy of a biregular automorphism by algebro-geometric methods. For instance, Example (\ref{tori1}) is immediate from this theorem; one may just compute $r_{p}(f)$ for $p=0$, $1$, $2$. Moreover, $d_p(f)$ and $r_p(f)$ can be regarded as {\it finer} invariants of $f$ than $h_{{\rm top}}(f)$, while geometric meanings become less apparent. 

\begin{corollary}\label{algebraicinteger}

(1) $h_{{\rm top}}(f) = 0$ if $\dim\, M = 1$, and also $h_{{\rm top}}(f) = 0$ 
for $f \in {\rm Aut}^0(M)$ (the identity component of ${\rm Aut}\, (M)$). For instance, $h_{{\rm top}}(f) = 0$ 
if $f \in {\rm Aut}\, ({\mathbf P}^d)$ or again if $f$ is a translation automorphism of 
a complex torus. 

(2) The topological entropy is the natural logarithm of an algebraic integer. 
\end{corollary}

Indeed, (1) is clear by Theorem (\ref{GromovYomdin}). Since the eigenvalues of $f^*\vert H^*(M, {\mathbf Z})$ are algebraic integers, (2) follows from Theorem (\ref{GromovYomdin}).

\begin{corollary}\label{concave}

(1) $d_0(f)=1$ and $d_l(f)=\, {\rm deg}\, f$, the topological degree of $f$. 

(2) The sequence $\{d_{p}(f)\}_{0 \le p \le l}$ is log-concave, i.e., 
$d_{p-1}(f)d_{p+1}(f) \leq d_p(f)^2$.

(3) $d_{p}(f) \ge 1$ for all $p$ and, 

(4) $h_{{\rm top}}(f) > 0$ if and only if $d_1(f) > 1$. 
\end{corollary}

(1) is clear by definition. We have $\delta_{p-1}(f)\delta_{p+1}(f) \le \delta_p(f)^2$ by the Hodge index theorem (when $M$ is projective and $\eta$ is chosen to be a Hodge metric) and by the Tessier-Khovanski inequality in general case. 
Then (2) follows from (1), and (2) implies (3), (4). 

Very brief outline of the proof of Theorem (\ref{GromovYomdin}) is as follows (\cite{Gr87}, \cite{DS05-1}). 
Note the {\it obvious} relation $(f^n)^* = (f^*)^n$ for 
$f \in {\rm Aut}\, (M)$. Then $d_p(f) = r_{p}(f)$ follows from linear algebra plus the Perron-Frobenious theorem on the linear maps preserving a strict convex cone. So, $d_{p}(f) = r_{p}(f) \le r^{{\rm even}}(f) \le r(f)$. The deepest part is $h_{{\rm top}}(f) \ge \log r(f)$ for any compact oriented Riemannian $C^{\infty}$-manifolds and surjective oriented $C^{\infty}$-map $f : M \to M$. This is due to Yomdin \cite{Yo87} (see also \cite{Gr87}). Gromov \cite{Gr03} 
proved the reverse inequality $h_{{\rm top}}(f) \leq {\rm lov}\, (f) = 
\log \max _{0\leq p\leq k} d_p(f)$. The essential part of this inequality is that if $F \subset M$ is $(n, \epsilon)$-separated, then the corresponding subset ${\mathbf F}^{(n)}$ in the graph 
$\Gamma_{f, n}$ is 
$(1, \epsilon)$-separated, and therefore the balls $\Gamma_{f, n} \cap B({\mathbf x}^{(n)}, \epsilon/2)$ (${\mathbf x}^{(n)} \in {\mathbf F}^{(n)}$) are mutually disjoint. This gives an obvious estimate of ${\rm Volume}\, (\Gamma_{f, n})$ from the below and leads the first inequality, via Lelong's theorem. ${\rm lov}\, (f) = {\rm max}_{0 \le p \le d} d_{p}(f)$ is non-trivial but doable by fairly straightforward computations of the differential forms, {\it at least when $f$ is holomorphic}. 

See also \cite{DHKK13} for derived categorical approach for the topological entropy. 

\subsection{Generalization for rational mappings after Dinh-Sibony} References of this subsection are \cite{Gu05}, \cite{DS05-1}, \cite{DS05-2} (see also \cite{CCG10}). 

Let $M$ be a compact K\"ahler manifold of dimension $l$, and $f: M \cdots\to M$ be a dominant rational selfmap. The topological entropy $h_{{\rm top}}(f)$ is defined in the same way as in the holomorphic case, just by considering the well-defined orbits, i.e., $\{f^k(x)\}_{k=1}^{n-1}$ with $f^k(x) \not\in I(f)$ at each $n$-th step. The logarithmic volume ${\rm lov}\, (f)$ is defined again in the same way as above, just by taking the graph $\Gamma_{f, n}^0$ over $M \setminus \cup_{k=0}^{n-1} I(f^k)$ at each $n$-th step (\cite{Gu05}, \cite{DS05-1}). 

The pullback operation $f^* : H^{p, p}(M) \to H^{p, p}(M)$ is well-defined if one uses {\it currents}. Let $\tilde{M}$ be a resolution of the indeterminacy of $f$ and $p_i : \tilde{M} \to M$ ($i=1$, $2$) be the natural projections. 
Then for any closed $(p, p)$-form $\alpha$, we define $f^*(\alpha) = (p_1)_{*}p_2^*(\alpha)$, where $p_2^*$ is the natural pullback as forms and $(p_{1})_{*}$ is the natural pushforward as currents, i.e., $\langle (p_{1})_{*}(p_2^*(\alpha)), \beta \rangle_M := \langle p_2^*(\alpha), p_{1}^{*}(\beta)\rangle_{\tilde{M}}$. The action $f^*$naturally descends to the linear action on $H^{p, p}(M)$. So, the definitions of $\delta_p(f)$ and the $p$-th dynamical degree $d_p(f)$ make sense without any change. $d_p(f)$ does not depend on the choice of $\eta$ for the same reason as before, again once the existence of the limit is guaranteed. However, there is one crucial difference from the holomorphic case; {\it $(f \circ g)^* \not= g^* \circ f^*$ and $(f^n)^* \not= (f^*)^n$ in general.} This already happens for the standard Cremona transformation ${\rm cr}_l$ of ${\mathbf P}^l$ and makes outlined proof in the holomorphic case delicate at all the places where we freely use them. For instance, 
in general, there is no way to compare $d_p(f)$ and $r_p(f)$. Dinh and Sibony (\cite{DS05-1}\, \cite{DS05-2}) proved:

\begin{theorem}\label{DinhSibony2}
Let $X$ be a compact K\"ahler manifold of dimension $l$ and $f : M \cdots\to M$ be a dominant, rational map ($=$ meromorphic map, by our convention). Then, Theorem (\ref{DinhSibony1}) holds for rational dominat maps and 
$$h_{{\rm top}}(f) \le {\rm lov}\, (f) = \log {\rm max}_{0 \le p \le d}\, d_{p}(f)\,\, .$$
Moreover $d_p(f)$ are birational invariants in the sense that $d_p(f) = d_p(\varphi \circ f \circ \varphi^{-1})$ for any birational map 
$\varphi : M \cdots\to M'$ between compact K\"ahler manifolds. 
\end{theorem}
On the other hand, $h_{{\rm top}}(f)$ is {\it not} a birational invariant by Guedj \cite{Gu05}: 
\begin{example}\label{Guedj}
Let $f : {\mathbf C}^2 \to {\mathbf C}^2$ be a morphism defined by 
$(x, y) \to (x^2, y+1)$. Then $f$ naturally extends to a rational selfmap $f_1$ of ${\mathbf P}^2$ and a holomorphic selfmap $f_2$ of ${\mathbf P}^1 \times {\mathbf P}^1$. Then, $h_{{\rm top}}(f_1) = 0$ but $h_{{\rm top}}(f_2) = \log 2 > 0$.
\end{example}
Because of the birational invariance of dynamical degrees, one can define dynamical degrees for a dominant selfmap of a singular compact complex space of class ${\mathcal C}$ in an obvious manner, but not the topological entropy 
in this way. For this reason, dynamical degrees fit well more with birational geometry than the entropy. They are also useful when we study {\it biregular} automorphisms in Problem (\ref{Explicit1}), as we shall see in concrete cases in Section 5. 

The essential part of the proof of Theorems (\ref{DinhSibony1}), (\ref{DinhSibony2}) and their variants later is a deep theory of semi-regularization of currents, very roughly, a method for approximating currents well by sequences of smooth forms. Once such semi-regularization results are well established, then the proof goes along the same line as in the holomorphic case {\it if one carefully replaces all necessary estimates for currents by those of semi-regularizing smooth forms}. 

For a rational dominant selfmap $f$, Corollary (\ref{algebraicinteger}) (2) {\it is expected to be true but unknown}. For instance, (2) is true for $d_1(f)$ if $f \in {\rm PsAut}\, (M)$. This is because then $(f^n)^* = (f^*)^n$ on $H^2(M, {\mathbf Z})$, hence $d_1(f) = r_1(f)$ for the same reason as in the holomorphic case. Corollary (\ref{concave}) (1), (2), (3), being free from $h_{{\rm top}}(f)$, is clearly true, but (4) is {\it not true} as Example (\ref{Guedj}) shows. 

\subsection{Relative dynamical degrees after Dinh-Nguy\^en-Truong} References of this subsection are \cite{DN11} and \cite{DNT11}. Corollary (\ref{entropybasic})(3) or dynamical degrees of the product map $f \times f' : X \times X' \to X \times X'$ suggests a good notion of relative dynamical degrees with nice properties. If exists, then it will provide a useful numerical criterion for the primitivity of a selfmap, as we shall test in some concrete cases in Section 5. 

{\bf Setting I.} {\it Let $f: M\cdots \to M\,\, ,\,\, g: B\cdots \to B$ be dominant rational maps such that $\pi \circ f=g\circ \pi$. Here $\pi : M \to B$ is a surjective holomorphic map between compact K\"ahler manifolds $M$ and $B$ of dimensions $l$ and $b$ (necessarily $l \ge b$) with K\"ahler forms 
$\eta_M$ and $\eta_B$.}

In Setting I, we define {\it the relative dynamical degrees} $d_p(f|\pi )$ by
\begin{eqnarray*}
d_p(f|\pi ) :=\lim _{n\rightarrow\infty}(\int_M (f^n)^*(\eta_M^p) \wedge \pi ^*(\eta_B^{b})\wedge \eta_M^{l-p-b})^{1/n}\,\, ,\,\, 0\leq p\leq l-b\,\, .
\end{eqnarray*}
This definition is due to Dinh and Nguy\^en \cite{DN11}. If we take $\eta_B$ so that $\eta_B^b$ is the Poincar\'e dual of a point and {\it virtually identify} all the fibers $M_b$ and regard then $f \vert M_b : M_b \cdots\to M_{g(b)}$ as the {\it virtual selfmap} of $M_b$, then $d_p(f|\pi)$ is the same form as the $p$-th dynamical degree of the {\it virtual} $f \vert M_b$. Note also that $\pi ^*(\eta_B^{b})\wedge \eta_M^{l-p-b}$ is a form as $\pi$ is {\it holomorphic} and $(f^n)^*(\eta_M^p)$ is a current of proper bidegree. So the integration in the right hand side makes sense. This is the reason why we assume that $\pi$ is holomorphic. The existence of the limit is again non-trivial, but settled by \cite{DN11} and \cite{DNT11}. 
The following fundamental result is due to Dinh-Nguy\^en-Truong (\cite{DN11}, \cite{DNT11}):
\begin{theorem}\label{DTT} In Setting I, for all $0 \le p \le l$, 
$$d_p(f)= \max_{j}d_j(g)d_{p-j}(f|\pi)\,\, .$$ 
Here $j$ runs through all the integers for which the integrations defining $d_j(g)$ and $d_{p-j}(f\vert \pi)$ are meaningful, i.e., $j$ runs through ${\rm max}\, \{0, p - l +b\} \le j \le {\rm min}\, \{p, b\}$. Moreover, $\{d_p(f|\pi)\}_p$ satisfy $d_{p-1}(f \vert \pi)d_{p+1}(f \vert \pi) \le d_{p}(f \vert \pi)^2$ (the log-concavity) and they are birational invariants in an obvious sense, within Setting I.
\end{theorem}

That $\pi$ is holomorphic in Theorem (\ref{DTT}) appears slightly restrictive, compared with usual situations, and the most natural setting is probably the following: 

{\bf Setting II.} {\it $\pi' : M' \cdots\to B'$ is a dominant rational map from an $l$-dimensional compact complex variety $M'$ of class ${\mathcal C}$ to a compact complex variety $B'$ of dimension $b$, equivariant with rational dominant selfmaps $f'$ and $g'$ of $M'$ and $B'$.}

In Setting II, $B'$ is of class ${\mathcal C}$ by the original definition of class ${\mathcal C}$, and therefore, there is a birational morphism $\varphi : B \to B'$ from a compact K\"ahler manifold $B$ as well. Then resolving the indeterminacy of the rational map $\varphi^{-1} \circ \pi'$ from $M'$ to $B$, we obtain a {\it holomorphic} surjective morphism between compact K\"ahler {\it manifolds} $\pi : M \to B$. Moreover, $\pi$ is equivariant to the rational dominant selfmaps $f$ and $g$ of $M$ and $B$, naturally induced from $f'$ and $g'$. This is exactly Setting I in Theorem (\ref{DTT}). By the birational invariance of dynamical degrees, we have $d_p(f) = d_{p}(f')$, $d_p(g) = d_{p}(g')$. Moreover, by the birational invariance of the relative dynamical degrees in Setting I in Theorem (\ref{DTT}), we {\it can define} $d(f'\vert \pi') := d(f \vert \pi)$ 
which is independent of the choice of models $\pi : M \to B$. Then, the equation in Theorem (\ref{DTT}) is nothing but the equation $d_p(f')= \max_{j}d_j(g')d_{p-j}(f'|\pi')$ in Setting II. Therefore: 
\begin{corollary}\label{generalDTT} Theorem (\ref{DTT}) is true also in 
the Setting II.
\end{corollary}

Note that $d_0(f'|\pi') =1$ 
and $d_{l-b}(f'|\pi' )$ is the topological degree of $f'\vert M_t : M_t \cdots\to M_{g'(t)}$ for a generic fiber $M_t'$ ($t \in B$). The log-concavity 
then implies that $d_p(f'|\pi') \ge 1$ for any meaningful $p$ as before.

Since only $d_0(f|\pi) = 1$ is the meaningful relative dynamical degree for a generically finite map, we obtain the following (\cite{DN11}, \cite{DNT11}) from Theorem (\ref{generalDTT}):
\begin{corollary}\label{genericallyfinite} 

(1) The dynamical degrees are invariant under any 
equivariant generically finite 
dominant maps, i.e., if $\pi : M \cdots\to  B$ is a generically finite dominant rational map equivariant to the rational dominant selfmaps $f$, $g$ of $M$, $B$, then $d_p(f)= d_p(g)$ for every $p$. 

(2) The topological entropy of dominant holomorphic selfmaps of compact K\"ahler manifolds are invariant under equivariant generically finite dominant rational maps. More precisely, in (1), if both $M$ and $B$ are compact K\"ahler manifolds and $f$ and $g$ are holomorphic, then $h_{{\rm top}}(f)=h_{{\rm top}}(g)$. 
\end{corollary} 

(2) follows from (1) and Theorem (\ref{GromovYomdin}). 

Our primary interest in Theorem (\ref{DTT}) is its applicability for primitivity of rational selfmaps. When $l = \dim\, M \le 3$, we can deduce the following fairly useful {\it numerical} criterion for the primitivity of $f \in {\rm Bir}\, M$ from Theorem (\ref{DTT}). (1) is known before Theorem (\ref{DTT}) and 
(2) is due to Truong and myself \cite{OT13-2}:

\begin{corollary}\label{TO1} Let $M$ be a compact K\"ahler manifold and $f \in {\rm Bir}\, (M)$. 

(1) Assume that $\dim\, M = 2$. Then $f$ is primitive if $d_1(f) > 1$. In particular, 
$f \in {\rm Aut}\, (M)$ is primitive if $h_{{\rm top}} (f) > 0$. 

(2) Assume that $\dim\, M = 3$. Then $f$ is primitive if $d_1(f) \not= d_2(f)$. 

\end{corollary}
 
Outline of (2) is as follows. Assume that $f$ is imprimitive. Then there are a compact K\"ahler manifold $B$, dominant rational maps $\pi : M \cdots \to B$, $g : B \cdots \to B$ such that $\pi \circ f=g\circ \pi$. Here $0 < \dim\, B < 3 =\dim\, M$. We consider the case $\dim\, B = 2$ (the case $\dim\, B = 1$ is similar). Then 
by Corollary (\ref{generalDTT}), we have 
$$d_1(f) = \max \{d_1(g), d_1(f|\pi) \}\,\, ,\,\, d_2(f)=\max \{d_1(f|\pi)d_1(g), d_2(g)\}\,\, .$$
Since $f$ and $g$ are birational, $d_3(f) = d_2(g) = 1$. 
Thus, by Corollary (\ref{generalDTT}), $1 = d_3(f) = d_2(g)d_1(f|\pi )$, hence, $d_1(f|\pi) = 1$. So, $d_1(f) = \max \{d_1(g), 1\} = d_2(f)$. 

\section{Surface automorphisms in the view of entropy} In this section, 
we take a closer look at surface automorphisms of positive entropy. They are primitive by Corollary (\ref{TO1})(1).  
We assume some familiarity with classification of surfaces. 
A good reference is \cite{BHPV04} with \cite{Do11} for rational surfaces. 
Throughout this section, $S$ is a smooth compact K\"ahler surface. 

\subsection{Surface automorphisms of positive entropy} 
We note that a birational automorphism $f \in {\rm Bir}\, (S)$ naturally induces a {\it biregular} automorphism of the minimal model $S_{{\rm min}}$ of the same dynamical degrees (Theorem (\ref{DinhSibony2})) if $\kappa\, (S) \ge 0$. In this way, one can almost recover from Theorem (\ref{deqizhang}) the following breakthrough observation due to Cantat \cite{Ca99}:

\begin{theorem}\label{cantat2} Assume that $S$ admits an automorphism $f \in {\rm Aut}\, (S)$ of positive entropy, i.e., $d_1(f) > 1$. Then $S$ is birational to either (i) ${\mathbf P}^2$, 
(ii) a K3 surface, (iii) a $2$-dimensional complex torus, or (iv) an Enriques surface. In the case (i), $S$ is a blow up of ${\mathbf P}^2$ at $10$ or more points, possibly infinitely near (\cite{Na60}). 
\end{theorem} 

Recall that $d_1(f) = r_1(f) >1$ is an algebraic integer (Corollary (\ref{algebraicinteger}) (2)). It turns out to be a special algebraic integer of 
even degree, called a {\it Salem number}: 
\begin{definition}\label{salem} An irreducible monic polynomial $S(x) \in {\mathbf Z}[x]$ is called a {\it Salem polynomial} if the complex roots are of the following form ({\it possibly} $d=1$):
$$a \in (1, \infty)\,\, ,\,\, 1/a \in (0, 1)\,\, ,\,\, \alpha_i, \overline{\alpha_i} \in S^1 := \{z \in {\mathbf C}\, \vert\, \vert z \vert = 1\} \setminus \{\pm 1 \} (1 \le i \le d-1)\,\, .$$ 
The unique root $a >1$ is called a {\it Salem number} of degree $2d$ ($= {\rm deg}\, S(x)$).
\end{definition}
The smallest {\it known} Salem number is the {\it Lehmer number} which is the unique root $>1$ of the following Salem polynomial of degree $10$: 
$$x^{10} + x^9 -x^7 -x^6 -x^5 -x^4 -x^3 + x + 1\,\, .$$ 
It is approximately $1.17628$ and conjectured to be the minimum among all Salem numbers. So far, this conjecture is neither proved nor disproved. 

The following theorem is due to McMullen (\cite{Mc02-1}, \cite{Mc02-2}, \cite{Mc07}):

\begin{theorem}\label{mcmullen} Let $f \in {\rm Aut}\, (S)$ and assume that $d_1(f) > 1$. Then, $d_1(f)$ is a Salem number, and $d_1(f)$ is always greater than or equal to the Lehmer number. 
\end{theorem} 
See also \cite{Ue10}, \cite{Og10-2}, \cite{Mc11-1}, \cite{Mc11-2}, \cite{Re12} for relevant results. 

\subsection{Examples of surface automorphisms 
of positive entropy} 
There is a huge number of works concerning automorphisms of surfaces. 
Here among many examples, I present four examples which are smoothly connected to the topics in the next sections. 

{\bf Example 1 - Rational surface automorphisms.} Recall that any birational automorphism of ${\mathbf P}^2$ is expressed by two rational functions of the affine coordinates $(x = x_2/x_1, y = x_3/x_1)$ of ${\mathbf P}_{[x_1 : x_2 : x_3]}^2$. Consider the birational automorphism of the following special form (\cite{Mc07}, also compare with an alternative form in \cite{BK06}):
$$f^*(x, y) := f_{(a, b)}^*(x, y) := (a + y, b + \frac{y}{x})\,\, ,\,\, a, b \in {\mathbf C}\,\, .$$
$I(f)$ is the set of coordinate points $\{e_1, e_2, e_3\}$ and $I(f^{-1})$ 
is $\{e_2, e_3, e_4\}$, where $e_4 := (a, b)$. Set $e_{k+4} := f^{k}(e_4)$. 
Choose $(a, b)$ so that $e_k \not\in e_1e_2 \cup e_2e_3 \cup e_3e_1$ for all 
$4 \le k \le n$ with $n \ge 10$ and $e_{n+1} = e_1$ (periodicity condition of indeterminacies). Then one can realize $f$ as an automorphism of $S = S_n$, the blowing-ups of ${\mathbf P}^2$ at the $n$ points $e_k$. This $f$ also realizes the Coxeter element $c_n$ of the Weyl group $W(E_n)$ in the sense that $f^* = c_n$ on $H^2(S, {\mathbf Z})$ under the natural identification $E_n = K_S^{\perp}$ and $W(E_n) < {\rm O}\,(H^2(S, {\mathbf Z}))$, hence of positive entropy. For instance, $f$ with $n = 10$ and $(a, b) = (0.4995..., -0.0837...)$ realizes the Lehmer first dynamical degree and $f$ with $n=11$ and $(a, b)$ approximately $(0.0444 -0.4422i, 0.0444 +0.4422i)$ has a Siegel domain (cf. Example 4). The actual construction in \cite{Mc07} is not merely the numerical one but is based on an explicit marked Torelli type result for log K3 surfaces $(S, C)$ with $C$ being a unique cuspidal rational curve in $\vert -K_S \vert$.  

{\bf Example 2 - Birational automorphisms after Diller-Favre.} Reference here is \cite{DF01}. Let $c \in {\mathbf C}$ and consider the birational automorphism of ${\mathbf P}^1 \times {\mathbf P}^1$ defined by the following affine form:
$$f_{c}(x, y) := (y+1 - c,  x \frac{y - c}{y +1})\,\, .$$ 
\cite{DF01} computes the first dynamical degree of $f_{c}$ and observes many interesting features, depending on $c \in {\mathbf C}$. For instance, if $c$ is irrational, then $d_1(f_{c}) = \eta$, the golden number. Note that the golden number is {\it not} a Salem number, so that $f_{c}$ with irrational $c$ can {\it never} be realized as a biregular automorphism of any smooth birational models.

{\bf Example 3 - Cayley's K3 surface after Festi, Garbagnati, van Geemen and Luijk.} In the long history of automorphisms of K3 surfaces or more specifically those of smooth quartic surfaces, Cayley seems the first who suggested the existence of automorphisms of infinite order. Here we explain his beautiful, very explicit construction, following a modern elegant account \cite{FGGL12}. This example will be also used to construct higher dimensional HK example in Section 6. 

Let $a_{ijk}$ ($1 \le i, j, k \le 4$) be $4^3$ generic complex numbers. Let us consider the following determinantal quartic surface in ${\mathbf P}^3$ with homogeneous coordinates ${\mathbf x} = [x_1 : x_2 : x_3 : x_4]$:
$$S_0 := ({\rm det}\, M_0 ({\mathbf x}) = 0) \subset {\mathbf P}_{\mathbf x}^3\,\, ,$$
where $M_0 = M_0 ({\mathbf x}) := (\sum_i a_{ijk}x_i)_{k, j}$ is the $4 \times 4$ matrix whose $(k,j)$ entry is $\sum_i a_{ijk}x_i$. By our genericity assumption, 
${\rm rank}\, M_0 ({\mathbf x}) = 3$ for all ${\mathbf x} \in S_0$ and $S_0$ 
is a smooth quartic K3 surface. One can also construct two more smooth determinantal quartic K3 surfaces from $a_{ijk}$:
$$S_1 := ({\rm det}\, M_1 ({\mathbf y}) = 0) \subset {\mathbf P}_{\mathbf y}^3\, \, ,\,\, S_2 := ({\rm det}\, M_2 ({\mathbf z}) = 0) \subset {\mathbf P}_{\mathbf z}^3\,\, .$$
Here $M_1 = M_1 ({\mathbf y}) := (\sum_i a_{ijk}y_j)_{i, k}$, $M_2 = M_2 ({\mathbf z}) := (\sum_i a_{ijk}z_k)_{j, i}$. 

Let $P_i$ be the cofactor matrix of $M_i$. Then, 
$$P_iM_i = M_iP_i = {\rm det}\,( M_i) \cdot I_4\,\, {\rm and}\,\, 
{\mathbf x}M_1({\mathbf y}) = (\sum_{i, j}a_{ijk}x_iy_j)_k = {\mathbf y}M_0({\mathbf x})^t\,\, .$$
Recall that ${\rm rank}\, (M_0({\mathbf x})) \ge 3$ for each ${\mathbf x} \in {\mathbf P}^3$ and the same for $M_1({\mathbf y})$, $M_2({\mathbf z})$. Thus, the $j$-th column $(p_{ij}({\mathbf x}))_{i}$ of $P_0 = P_0({\mathbf x})$ gives a {\it Cremona transformation} $\varphi_0 : {\mathbf P}_{{\mathbf x}}^3 \cdots\to {\mathbf P}_{{\mathbf y}}^3$ that maps $S_0$ to $S_1$, hence, an isomorphism $\varphi_0\vert S_0 : S_0 \to S_1$. In the same way, we have two more Cremona transformations $\varphi_1 : {\mathbf P}_{{\mathbf y}}^3 \cdots\to {\mathbf P}_{{\mathbf z}}^3$, $\varphi_2 : {\mathbf P}_{{\mathbf z}}^3 \cdots\to {\mathbf P}_{{\mathbf x}}^3$ and isomorphisms $\varphi_1 \vert S_1 : S_1 \to S_2$, 
$\varphi_2\vert S_2 : S_2 \to S_0$. In this way, 
we obtain an {\it explicit} automorphism of $S_0$: $g := \varphi_2 \circ \varphi_1 \circ \varphi_0$. This is the automorphism that {\it Cayley first found around 1870 and said that "The process may be indefinitely repeated"} (\cite{Ca18}, \cite{FGGL12}). 

As observed by \cite{FGGL12}, in modern terminologies, a characterization of linear determinantal varieties (\cite{Be00}) with our genericity assumption says that $S_0$ is nothing but a K3 surface with 
$${\rm NS}\, (S_0) = ({\mathbf Z}h_1 \oplus {\mathbf Z}h_2, \left(\begin{array}{rr} 
4 & 6\\
6 & 4\\
\end{array} \right) ) = ({\mathbf Z}[\eta], 4{\rm Nm}(*))\,\, ,\,\, {\rm Nm}(a + b\eta) = a^2 + ab - b^2\,\, .$$
Here $\eta$ is the golden number and the lattice identification is given by $h_1 \leftrightarrow 1$ 
and $h_2 \leftrightarrow \eta^2$. 
Under this identification, the action of $g$ on ${\rm NS}\, (S)$ is the multiplication by $\eta^6$ on ${\mathbf Z}[\eta]$. So, as predicted by Cayley, $g$ is actually of infinite order and $h_{{\rm top}}(g) = \log \eta^6 > 0$. \cite{FGGL12} further shows that ${\rm Aut}\, (S_0) = \langle g \rangle$. They also give the explicit {\it integers} $a_{ijk} \in {\mathbf Z}$ with desired properties. I re-discovered Cayley's automorphism in answering a question of Kawaguchi (\cite{Og12-1}, see also \cite{BHK13}):
  
\begin{theorem}\label{free} Let $W$ be a smooth compact K\"ahler surface with automorphism $f$ such that $f$ is of positive entropy and has no fixed point. Then $W$ is birational to a projective K3 surface, and the pair of Cayley's K3 surface $S_0$ and its automorphism $g$ is one of such examples.  
\end{theorem}

{\bf Example 4 - Non-projective K3 surface automorphism with Siegel domain after McMullen.} Let $f$ be an automorphism of a smooth surface $S$. We call a domain 
$U \subset S$ a 
{\it Siegel domain} of $f$ if $f(U) = U$ and $U$ is biholomorphic to the $2$-dimensional unit disk $\Delta^2$ with coordinates $(z_1, z_2)$ such that the induced action of $f$ on $\Delta^2$ is of the form $f^*(z_1, z_2) = (\alpha_1z_1, \alpha_2z_2)$ for some multiplicatively independent complex numbers $\alpha_1$ and $\alpha_2$ on the unit circle $S^1$, i.e., $\alpha_1^{m_1}\alpha_2^{m_2} \not= 1$ for any integers $(m_1, m_2) \not= (0, 0)$ and $\vert \alpha_1 \vert = \vert \alpha_2 \vert = 1$. 

If $S$ is a K3 surface and $f$ is an automorphism with Siegel domain as above, then $f^*\sigma_S = \alpha_1\alpha_2 \sigma_S$. Here $\sigma_S \not= 0$ is a global holomorphic $2$-form on $S$. Note that $\alpha_1\alpha_2$ is 
not root of unity. Thus $S$ is necessarily non-projective, as the pluri-canonical representation of ${\rm Bir}\, (M)^* \vert H^0(M, {\mathcal O}_M(mK_M))$ is always finite if $M$ is projective (\cite{Ue75}). The next very surprising result due to McMullen (\cite{Mc02-2}) gave me a strong motivation 
to study birational automorphisms from the view of this note: 
\begin{theorem}\label{mcmullen3} 
There is a K3 surface $S$ of Picard number $0$ with ${\rm Aut}\, (S) = \langle f \rangle$ such that $f$ is of positive entropy and has a Siegel domain. In particular, the canonical representation $f^* \vert H^0({\mathcal O}(K_S)) = f^* \vert H^0(\Omega_S^2)$ is of infinite order, and there is no point $Q \in S$ such that the orbit ${\rm Aut}\, (S) \cdot Q$ is topological dense (even though $f$ is of positive entropy). Slightly more explicitly, one of such $(S, f)$ is realizable so that the characteristic polynomial of $f^* \vert H^2(S, {\mathbf Z})$ 
is the following Salem polynomial of degree $22$; 
$$S_{22}(X) := x^{22} + x^{21} -x^{19} - 2x^{18} -3x^{17} -3x^{16} -2x^{15} + 2x^{13} + 4x^{12}$$
$$+ 5x^{11} + 4x^{10} + 2x^{9} -2x^{7} -3x^{6} -3x^{5} -2x^{4} -x^{3} + x + 1\,\, .$$
In this case, $h_{{\rm top}}(f) = \log a$, where $a$ is the Salem number of $S_{22}(x)$, approximately, $1.37289$. 
\end{theorem}

Unlike the examples above, construction is highly implicit, based on the surjectivity of the period map and golobal Torelli theorem for K3 surfaces, and the existence of Siegel domain is based on a deep transcendental number theoretical result (\cite{Mc02-2}), from which one can deduce the transcendency of $\pi$ and $e$ in one line.

We close this section with a few remarks relevant to Theorem (\ref{mcmullen3}):

\begin{remark}\label{ratsiegel} As mentioned, there are smooth rational surfaces with an automorphism with Siegel domain (\cite{Mc07}, \cite{BK09}). Rational surfaces are always projective, but this does not contradict the finiteness of pluri-canonical representation, because $\kappa (S) = -\infty$, i.e., $H^0(S, {\mathcal O}_S(mK_S)) = 0$ for all $m > 0$. 
\end{remark}

\begin{remark}\label{voisin} Let $S$, $f$ be as in Theorem 
(\ref{mcmullen3}) and $P$ be the center of the Siegel domain.
Let $M$ be the blowing-ups of $N := S \times S$, first at 
the intersection point 
$(P, P)$ of 
$S \times \{P\}$, $\{P\} \times S$, 
the diagonal $\Delta$, and the graph $\Gamma_f$, and next along the proper transforms of $S \times \{P\}$, $\{P\} \times S$, $\Delta$, $\Gamma_f$. Then $M$ is a simply-connected compact K\"ahler fourfold which can not be deformed into projective manifolds under any small proper deformation (\cite{Og10-1}, also compare with \cite{Vo04}). 
\end{remark} 

\begin{remark}\label{picard} Let $S$ be a projective K3 surface. Then 
$1 \le \rho(S) \le 20$, and projective K3 surfaces with $\rho(S) \ge \rho$ form countable union of $(20 -\rho)$-dimensional families. The automorphism group of $S$ {\it tends} to be larger if $\rho(S)$ becomes larger (see \cite{Og03} for the precise statement in terms of deformation). If $\rho(S) = 20$, for instance if $S$ is the Fermat quartic surface, then ${\rm Aut}\, (S)$ contains the free subgroup 
${\mathbf Z} * {\mathbf Z}$ with many elements of positive entropy, and the orbit ${\rm Aut}\, (S) \cdot P$ is topologically dense in $S$ for generic $P \in S$ (\cite{Og07-2}, \cite{Ca01} see also \cite{SI77}). 
\end{remark}

\section{Rational and CY threefolds with primitive automorphisms of positive entropy}

\subsection{Biregular automorphisms vs. birational automorphisms} 
Some experiences show that in dimension $\ge 3$, the biregular automorphisms tend to be drastically fewer than birational automorphisms. So, finding manifolds with ``interesting" biregular automorphisms is more challenging in some sense. Here I present a few examples of this tendency.

{\bf Example 1 - CY manifolds in Fano manifolds.} Cayley's K3 surfaces are smooth anti-canonical members of the Fano threefold ${\mathbf P}^3$. Smooth anti-canonical members of higher dimensional Fano manifolds are CY manifolds. However, 

\begin{theorem}\label{kollar} Let $l \ge 3$ and $M$ be a smooth member of $\vert -K_V \vert$ of a smooth Fano manifold $V$ of dimension $l+1$. Then $M$ is a CY manifold of dimension $l \ge 3$, but $\vert {\rm Aut}\, (X) \vert$ is finite.
\end{theorem}

Lefschetz hyperplane theorem shows that $\iota^* : H^2(V, {\mathbf Z}) \simeq H^2(M, {\mathbf Z})$ under the inclusion $\iota : M \to V$ for $l \ge 3$. Koll\'ar (\cite{Bo91}) shows that $\overline{{\rm Amp}\, (M)} \simeq \overline{{\rm Amp}\, (V)}$ under $\iota^*$. Hence $\overline{{\rm Amp}\, (M)}$ is a finite rational polyhedral cone as so is $\overline{{\rm Amp}\, (V)}$. This implies $\vert {\rm Aut}\, (M) \vert < \infty$ by Remark (\ref{ample}). 

{\bf Example 2 - CY manifolds of smaller Picard numbers.} Recall that Cayley's K3 surfaces are of Picard number $2$ and have automorphism of positive entropy. On the other hand, we have (\cite{Og12-2}; see also \cite{LP12}, \cite{LOP13}):
\begin{theorem}\label{oguiso12} $\vert {\rm Aut}\, (M) \vert < \infty$ for an odd dimensional CY manifold of $\rho(M) = 2$. 
\end{theorem}

\begin{example}\label{oguisoex}
Let $M$ be a general complete intersection in ${\mathbf P}^3 \times {\mathbf P}^3$ of $2$ hypersurfaces of bidegree $(1,1)$ and a hypersurface of bidegree $(2,2)$. Then $M$ is a CY threefold of Picard number $2$ (hence an example of both Theorems (\ref{kollar}), (\ref{oguiso12})). Let $\iota_k$ be the covering involution of the $k$-th projection ${\rm pr}_k : M \cdots\to {\mathbf P}^3$ of degree $2$. Set $f := \iota_2 \circ \iota_1 \in {\rm Bir}\, (M)$. Then, $d_1 (f) = 17 + 12\sqrt{2} >1$, $\langle f \rangle \simeq {\mathbf Z}$ and $[{\rm Bir}\, (M) : \langle f \rangle ] < \infty$ even though ${\rm Aut}\, (M)$ is finite. 
\end{example}

Recall that ${\rm Bir}\, (M) = {\rm PsAut}\, (M)$ for CY manifolds. 
Example (\ref{oguisoex}) is an application of Theorem (\ref{kawamata1}) with an explicit analysis of the movable cone (\cite{Og12-2}). 

The following theorem also shows a sharp contrast in dimension $2$ and 
$\ge 3$:
\begin{theorem}\label{oguisocantat} Let $l \ge 2$ and $M = (2, \cdots ,2) \subset ({\mathbf P}^{1})^{l +1} = {\mathbf P}_1^1 \times {\mathbf P}_2^1 \times \cdots \times {\mathbf P}_{l+1}^{1}$ be a smooth generic hypersurface of multi-degree $(2, \cdots ,2)$. Then $M$ is a K3 surface if $l = 2$ and a CY manifold of dimension $l$ if $l \ge 3$, and:

(1) If $l =2$, then ${\rm Bir}\, (M) = {\rm Aut}\,(M) = \langle \iota_1, \iota_2, \iota_{3} \rangle \simeq {\mathbf Z}_2 * {\mathbf Z}_2 * {\mathbf Z}_2$, while

(2) If $l \ge 3$, then 
${\rm Aut}\,(M) = \{id_M \}$ and  
${\rm Bir}\, (M) = \langle \iota_1, \cdots, \iota_{l+1} \rangle 
\simeq {\mathbf Z}_2 * \cdots * {\mathbf Z}_2$ ($(l+1)$-times free product). 

Here $\iota_k$ is the covering involutions of the natural projection to the product $({\mathbf P}^1)^l$ in which the $k$-th factor ${\mathbf P}_k^1$ of $({\mathbf P}^1)^{l+1}$ removed. Moreover, there are (many) elements $f$ with $d_1(f) >1$.
\end{theorem}

So, ${\rm Bir}\, (M)$ becomes larger and larger according to the dimension, but ${\rm Aut}\, (M)$ suddenly disappears in dimension $\ge 3$. This is proved by Cantat and myself \cite{CO11}. The essential algebro-geometric part of (2) is that in $l \ge 3$, the covering involutions $\iota_k$ are (no longer automorphism but) a birational involutions and at the same time all possible flops of $M$. We then apply Theorem (\ref{kawamata1}).  

{\bf Example 4 - CY manifold automorphisms of positive entropy.} Only one "series" of examples with automorphisms of positive entropy that I know is:

\begin{theorem}\label{oguisocantat2} 
Let $M$ be the universal cover of the punctual Hilbert scheme 
${\rm Hilb}^l\, (S)$ of length $l \ge 2$ of an Enriques surface 
$S$. Then $M$ is a CY manifold of dimension $2l$, and $M$ admits (many) biregular automorphisms of positive entropy if $S$ is generic.
\end{theorem} 
See \cite{OS11} and \cite{CO11}. It is interesting to ask:
\begin{question}\label{enriques}
Does a CY manifold $M$ in Theorem (\ref{oguisocantat2}) admit a primitive automorphism of positive entropy?
\end{question}

{\bf Higher dimensional rational manifold automorphisms of positive entropy.} Finding ``interesting" {\it biregular} automorphisms of rational manifolds seems much more difficult. Surprisingly, the following most basic question, posed by Bedford, is still unsolved (See also \cite{Tr12} for many negative evidences): 

\begin{question}\label{bedford} 
Is there a {\it biregular} automorphism of positive entropy on a smooth rational threefold obtained by blowing-ups of 
${\mathbf P}^3$ along smooth centers? 
\end{question}

On the other hand, in \cite{PZ11}, \cite{BK14}, \cite{BCK13}, \cite{BDK14}, there are constructed many examples of rational manifolds with interesting pseudo-automorphisms, by generalizing constructions of rational surface automorphisms. Especially, the following result due to Bedford-Cantat-Kim \cite{BCK13} is quite remarkable and also strongly supports a negative answer for Question (\ref{bedford}):

\begin{theorem}\label{PrimitiveNotBiregular} 
There is a smooth rational threefold $M$, obtained by blowing-ups of 
${\mathbf P}^3$ at smooth centers, with $f \in {\rm PsAut}\, (M)$ such that $f$ is primitive, 
of $d_1(f) > 1$, and $(M, f)$ has no equivariant smooth birational model $(M', f')$ with $f' \in {\rm Aut}\, (M')$. 
\end{theorem}

Their $f$ does not preserve even any foliation. In their construction, $\pi : M \to {\mathbf P}^3$ is a partial resolution of indeterminacy of 
$\tau \circ {\rm cr}_3 \in {\rm Bir}\, ({\mathbf P}^3)$ and their iterates, for suitably chosen $\tau \in {\rm PGL}\, (4, {\mathbf C})$ with periodicity conditions of indeterminacy. Then the birational automorphism $f =  \pi \circ (\tau \circ {\rm cr}_3) \circ \pi^{-1}$ becomes a pseudo-automorphism of $M$. In their construction, there is a rational surface $S \subset M$ preserved by $f^4$ in the sense that $f^4 \vert S \in {\rm Bir}\, (S)$. Their crucial observation for the non-existence of $(M', f')$ is that $d_1(f^4 \vert S)$ is {\it neither} a Salem number {\it nor} $1$. On the other hand, if $f$ could be regularized, then so is $f^4 \vert S$, possibly on other regularized models on which $S$ survives. 
But, then $d_1(f^4 \vert S)$ must be a Salem number or $1$ by the birational invariance of the dynamical degrees (Theorem (\ref{DinhSibony2})) and by Theorem (\ref{mcmullen}). 

\subsection{First examples of rational and CY threefolds with primitive automorphisms of positive entropy} Main reference is \cite{OT13-2}. The essential idea is the quotient construction from a manifold $M$ with rich automorphisms: If $G < {\rm Aut}\, (M)$ is a ``small" finite subgroup with ``big" normalizer 
$N < {\rm Aut}\, (M)$, then the ``big" group $N/G$ acts biregularly on the quotient variety $M/G$ and on its equivariant resolution as well.

Our actual construction is as follows. Let $E_{\tau} = {\mathbf C}/({\mathbf Z} + {\mathbf Z}\tau)$ be the elliptic curve of period $\tau$. There are exactly two elliptic curves with 
a Lie automorphism other than $\pm 1$. They are $E_{\sqrt{-1}}$ and $E_{\omega}$, where $\omega := (-1 + \sqrt{-3})/2$. Let $X_4$ (resp. $X_6$, resp. $X_3$) be the canonical resolutions of the quotient threefolds 
$$E_{\sqrt{-1}} \times E_{\sqrt{-1}} \times E_{\sqrt{-1}} /\langle 
\sqrt{-1} I_3 \rangle \,\, , 
\,\, E_{\omega} \times E_{\omega} \times E_{\omega}/ \langle -\omega I_3 \rangle\,\, , E_{\omega} \times E_{\omega} \times E_{\omega}/ \langle \omega I_3 \rangle$$
i.e., the blow up at the maximal ideals of singular points. As is well known, $X_3$ is a CY threefold (\cite{Be83}). It is analytically rigid, but plays 
an important role in the classification of CY threefolds in the view of the second Chern class (\cite{Og93}, \cite{OS01}). Our $X_3$, $X_6$, $X_4$ provide the {\it first} examples of a Calabi-Yau threefold and smooth rational threefolds with {\it primitive biregular} automorphisms of {\it positive entropy}:

\begin{theorem}\label{firstexample} 

(1) Both $X_6$ and $X_4$ are rational. 

(2) Moreover, $X_3$, $X_6$, $X_4$ admit {\it primitive} biregular automorphisms of {\it positive entropy}. 
\end{theorem}

(1) is proved by Truong and myself for $X_6$ (\cite{OT13-2}), and by Colliot-Th\'el\`ene (\cite{CTh13}) for $X_4$ via \cite{COT13} both answering a question of Ueno and Campana (\cite{Ue75}, \cite{Ca12}). (2) is proved by Truong and myself (\cite{OT13-2}). 

The most crucial parts are the rationality of $X_6$ and $X_4$ and finding primitive automorphisms. One of the key steps for the rationality is the following result (\cite{OT13-2}, \cite{COT13}) shown via determination of the rational function fields:

\begin{theorem}\label{AbstractField} 
Let $(s, t, z, w)$ be the standard affine coordinates of ${\mathbf C}^4$. Then:

(1) $X_4$ is birational to the hypersurface $H_4$ in ${\mathbf C}^4$ defined by
$$(t^2 - z)(s^2- w^3) = (s^2 - w)(t^2 - z^3)\,\, .$$

(2) $X_6$ is birational to the hypersurface $H_6$ in ${\mathbf C}^4$ defined by
$$(w^3 - 1)(t^2 -1) = (z^3 - 1)(s^2 -1)\,\, .$$
\end{theorem}
The projection $p_{34} : (t,s, z, w) \mapsto (z, w)$ gives the conic bundle 
structures on $H_4$ and $H_6$;  
$p_{34} : H_4 \rightarrow {\mathbf C}^2$, $p_{34} : H_6 \rightarrow {\mathbf C}^2$. It is clear that $(t,s, z, w) = (1,1, z, w)$ is a section of $p_{34} : H_6 \rightarrow {\mathbf C}^2$, and therefore $H_6$ is rational. $p_{34} : H_4 \rightarrow {\mathbf C}^2$ does not admit a rational section. However, 
Colliot-Th\'el\`ene \cite{CTh13} shows that the conic bundle $p_{34} : H_4 \rightarrow {\mathbf C}^2$ is birational to the conic bundle $p_{34} : (H_4)' \rightarrow {\mathbf C}^2$ over the same base. Here $(H_4)'$ is the affine hypersurface defined by $t^2 - zs^2 -w = 0$. This process is not explicit, but a consequence of the fact that these two conic bundles define the same element of the Brauer group ${\rm Br}\, ({\mathbf C}(z, w))$ of the base space ${\mathbf C}^2$ (\cite{CTh13}). $(H_4)'$ is rational as $w = t^2 -zs^2$, whence so is $H_4$.

Let us give an example of primitive biregular automorphisms of positive entropy of $X_3$, $X_6$, $X_4$. Let us consider the matrix
$$P = P_a = \left(\begin{array}{rrr}
0 & 1 & 0\\
0 & 0  & 1\\
-1 & 3a^2 & 0\\
\end{array} \right)\,\, ,$$
where $a$ is an integer such that $a \ge 2$. Since ${\det}\, P = 1$, $P$ naturally defines Lie automorphisms $g_3 = g_6$, $g_4$ of $(E_{\omega})^3$ and $(E_{\sqrt{-1}})^3$.  
By the construction and by the universality of blowing-up, the automorphisms $g_k$ descends to the {\it biregular} automorphisms, say $f_k$, of $X_k$. The eigenvalues of $\alpha$, $\beta$, $\gamma$ of $P$ are real numbers with $\vert \alpha \vert > \vert \gamma \vert > 1 > \vert \beta \vert$. By using $\alpha$, $\beta$, $\gamma$, one can compute that $d_1(g_k) = \alpha^2$ and $d_2(g_k) = \alpha^2\gamma^2$. Thus $d_2(g_k) > d_1(g_k) > 1$. Then 
$d_2(f_k) > d_1(f_k) > 1$ by Corollary (\ref{genericallyfinite}). Hence $f_k$ provide desired automorphisms by Corollary (\ref{TO1})(2). 

\begin{question}\label{geom} It is interesting to connect $X_6$ and $X_4$ to ${\mathbf P}^3$ by explicit blowing-ups and blowing-downs along smooth centers. This is in principle possible by \cite{AKMW02}. In the view of Question (\ref{bedford}), it is quite interesting to see if one can obtain $X_6$ and/or $X_4$ only by blowing-ups of ${\mathbf P}^3$ along smooth centers or not. 
\end{question}

\section{Birational automorphisms of HK manifolds}

\subsection{Some generalities}

We assume some familiarity with basics on HK manifolds. Excellent references are \cite{GHJ03}, Part III before Verbitsky's Torelli theorem (\cite{Ve09}) and \cite{Ma11} after that. We only recall that any HK manifold admits a non-degenerate integral symmetric bilinear form, called Beauville-Bogomolov form (BB-form, for short), $b_M(*, **) : H^2(M, {\mathbf Z}) \times H^2(M, {\mathbf Z}) \to {\mathbf Z}$ with signature $(3, b_2(M) -3)$, being compatible with Hodge decomposition and invariant under deformation and ${\rm Bir}\, (M)$. We denote by $S^{[n]} = {\rm Hilb}^n(S)$ the Hilbert scheme of the $0$-dimensional closed subschemes of lengths $n \ge 2$ on a K3 surface $S$. $S^{[n]}$ is a HK manifold of dimension $2n$ and of $\rho(S^{[n]}) = \rho(S) +1$. Through the Hilbert-Chow morphism $S^{[n]} \to {\rm Sym}^n\, S$, we have a natural identification as ${\mathbf Z}$-modules, $H^2(S^{[n]}, {\mathbf Z}) = H^2(S, {\mathbf Z}) \oplus {\mathbf Z}e$,  
where $e = [E]/2$, the half of the exceptional divisor $E$ of the Hilbert-Chow morphism. Under the BB-form $b(*, **)$, the above isomorphism is also an isometry with $(e^2) = -2(n-1)$.

\subsection{Rough structure theorem on birational automorphisms}

Let $M$ be a HK manifold of $\dim\, M = 2n$. We have the following Tits' alternatives:

\begin{theorem}\label{birhk} 
Let $G < {\rm Bir}\, (M)$. Then: 

(1) For each $(M, G)$, either one of the following two holds: 

(i) $G$ is an almost abelian group of rank $r$, i.e., $G$ is isomorphic to ${\mathbf Z}^r$ ($r \ge 0$) up to finite kernel and cokernel, or 

(ii) $G$ is essentially non-commutative, i.e., $G$ contains a free subgroup ${\mathbf Z} * {\mathbf Z}$.

(2) (ii) happens only if $M$ is projective and $\rho(M) \ge 3$. Moreover, in case (ii), there are (many) $f \in G$ such that $d_1(f) > 1$. In particular, if in addtion $G \subset {\rm Aut}\, (M)$, then there are (many) $f \in G$ with $h_{{\rm top}}(f) > 0$ in the case (ii). 
\end{theorem}

This is proved by \cite{Og07-1}, \cite{Og07-2}, \cite{Og08}. The essential part is as follows. The natural representation $r : G \to {\rm GL}\, (H^2(M, {\mathbf Z}))$ has a finite kernel (Remark (\ref{ample}) for projective case and \cite{Hu03} for general case). Thus $G$ is well approximated by the image $G^* := r(G)$. Then, the fundamental result of Tits (Theorem (\ref{tits})) reduces the problem to showing that if $G^*$ is virtually solvable, then it is an almost abelian group of finite rank. This can be done by using the additional strong condition 
that $G^*$ is a subgroup of ${\rm O}_{{\rm Hodge}}(H^2(M, {\mathbf Z}))$.
\begin{theorem}\label{tits} Any group $H < {\rm GL}\, (n, k)$ ($k$ is a field of characteristic $0$) satisfies either one of the following two:

(1) $H$ has a solvable subgroup of finite index (virtually solvable), or

(2) $H$ is essentially non-commutative, i.e., $H$ contains ${\mathbf Z} * {\mathbf Z}$.
\end{theorem}

\begin{remark}\label{titsgeneral} Tits' alternative type result with the same form as in Theorem (\ref{birhk}) does not hold in general. However, some meaningful different formulation is proposed and proved for the biregular automorphism group of any compact K\"ahler manifold (\cite{KOZ09}, \cite{Zh09-2}, see also \cite{DS04}, \cite{Di12}).
\end{remark}

One can also compute the dynamical degrees and entropy (\cite{Og09}):

\begin{theorem}\label{entropyhk} For any $f \in {\rm Aut}\, (M)$ 
of any HK manifold $M$, the dynamical degrees $d_k(f)$ are all Salem numbers or $1$. More precisely, $d_{2n -k}(f) = d_{k}(f) = d_1(f)^k$
for all $0 \le k \le n = \dim\, M/2$. In particular, if $d_1(f) >1$, then
$$1 = d_0(f) < d_1(f) < \cdots < d_{n-1}(f) < d_n(f) > d_{n+1}(f) > \cdots 
> d_{2n}(f) = 1\,\, ,$$
and $h_{{\rm top}}(f) = n\log\, d_1(f) >0$ (resp. $0$) if $d_1(f) >1$ (resp. $d_1(f) = 1)$.  
\end{theorem}

\subsection{A few examples}

By the definition of $S^{[n]}$, we have a natural inclusion 
${\rm Aut}\, (S) \subset {\rm Aut}\, (S^{[n]})$. This shows that if ${\rm Aut}\, (S)$ is infinite, then so is ${\rm Aut}\, (S^{[n]})$. So, contrary to the case of CY manifolds and rational manifolds, there are many examples of HK manifold with many biregular automorphisms. 

{\bf Example 1 - Non-projective primitive automorphism of positive entropy.}
Let $(S, f)$ be as in Theorem (\ref{mcmullen3}). 
Set $M = S^{[n]}$ and denote by $f_M \in {\rm Aut}\, (M)$ the automorphism naturally induced by $f$. We have $d_1(f_M) = d_1(f) = a$, the Salem number, and therefore $h_{{\rm top}}(f_M) = n\log\, a >0$. Since $S$ has no non-constant global meromorphic function (as $\rho(S) =0$), the same is true for $M$. Then $M$ has no rational fibration (\cite{COP10}). Hence $f_M$ is primitive as well. We also see that ${\rm Bir}\, (M) = {\rm Aut}\, (M) = \langle f_M \rangle \simeq {\rm Aut}\, (S) \simeq {\mathbf Z}$. 

{\bf Example 2 - The case where Picard number $2$.} We have the following:

\begin{theorem}\label{picard1}

(1) Let $S$ be a projective K3 surface with $\rho(S) = 1$. 
Then, $\rho(S^{[n]}) = 2$ but ${\rm Bir}\, (S^{[n]})$ is a finite group. 

(2) There is a projective HK fourfold $M$ deformation equivalent to $S^{[2]}$ 
such that $\rho(M) = 2$, ${\rm Aut}\, (M) = {\rm Bir}\, (M)$ is almost abelian group of rank $1$ with element of positive entropy. More specifically, $M$ with ${\rm NS}\, (M) \simeq ({\mathbf Z}[\eta], 4{\rm Nm}(*))$, "the same N\'eron-Severi lattice as Cayley's K3 surface", gives such an example. In particular, $2$ is the minimal Picard number of projective HK manfolds of dimension $\ge 4$ with automorphism of positive entoropy (cf. Remark (\ref{ample})). 
\end{theorem}

(1) is observed by \cite{Og12-2}. Unlike Cayley's K3 surfaces, our $M$ in (2) 
is highly non-constructible. {\it However, it is likely true that $M$ in (2) has a primitive automorphism of positive entropy (not yet settled).} 

{\bf Example 3 - Projective HK manifold of Picard number $3$.} Let $S \subset {\mathbf P}^3$ be a smooth quartic surface. Then for two general points $P, Q$ in $S$, the line $PQ$ in ${\mathbf P}^3$ meets $S$ in four points, say, $P$, $Q$, $P'$, $Q'$. The correspondence $\{P, Q\} \mapsto \{P', Q'\}$ defines a birational automorphism $\iota_{S}$ of $S^{[2]}$ of order $2$, called the {\it Beauville involution} (\cite{Be83}). If $S$ has no line, then $\iota_{S}$ is biregular. Note that $\iota_S \in {\rm Bir}\, (S^{[2]}) \setminus {\rm Aut}\, (S)$ under ${\rm Aut}\, (S) \subset {\rm Aut}\, (S^{[2]})$. 

Let $S$ be a Cayley's K3 surface. Identifying $S = S_{0} \subset {\mathbf P}^3$, our $S$ has three different embeddings $\Phi_k : S \to S_k \subset {\mathbf P}^3$ ($k = 0$, $1$, $2$) under the notation in Example 3 in Subsection (4.2). Let $\iota_k$ be the Beauville involution with respect to the embedding $\Phi_k$. 
We have the following theorem similar to Theorem (\ref{oguisocantat}):
\begin{theorem}\label{cayleyhk} Let $S$ be a Cayley's K3 surface. Then,
$${\rm Bir}\, (S^{[2]}) = {\rm Aut}\, (S^{[2]}) = \langle \iota_0, \iota_1, 
\iota_2 \rangle\,\, {\rm and}\,\, g = \iota_0 \circ \iota_1 \circ \iota_2\,\, ,$$ under the natural inclusion $\langle g \rangle = {\rm Aut}\, (S) \subset {\rm Aut}\, (S^{[2]})$.
Moreover, ${\rm Aut}\, (S^{[2]})$ has a subgroup isomorphic to the free product ${\mathbf Z} * {\mathbf Z}$, hence admits an automorphism of positive entoropy 
(Theorem (\ref{birhk})). In particular, $[{\rm Aut}\, (S^{[2]}) : {\rm Aut}\, (S) ] = \infty$ and $3$ is the minimal Picard number of projective HK manfolds of dimension $\ge 4$ with essentially non-commutative automorphism group.
\end{theorem}

One of interesting fact is that we have the second factorization $g = \iota_0 \circ \iota_1 \circ \iota_2$ in $S^{[2]}$, which looks 
similar to , but completely different from, the factorization that Cayley found in ${\mathbf P}^3$ (Example 3 in Subsection (4.2)). Another interesting fact is that ${\rm Aut}\, (S^{[2]})$ becomes much bigger than ${\rm Aut}\, (S)$ in this example, which makes a sharp contrast to the following open question, called the naturality question, posed by 
Boissi\`ere and Sarti (\cite{Bo11}, \cite{BNS11}):

\begin{question}\label{bs} Is ${\rm Aut}\, (S) = {\rm Aut}\, (S^{[m]})$ under the natural inclusion for $m \ge 3$?
\end{question} 

{\bf Acknowledgements.} First of all, I would like to express my sincere thanks to Professor Yujiro Kawamata for his continuous, warm encouragement, support and proper advices, both in mathematics and in life, since I was his graduate student on 1987. I would like to express my thanks to all my collaborators and teachers, especially Professors Tuyen Truong, De-Qi Zhang, Jun-Muk Hwang and Professors Fabrizio Catanese, Akira Fujiki, Heisuke Hironaka, JongHae Keum, Nessim Sibony, Tetsuji Shioda, Shing-Tung Yau, 
and Late Professors Eiji Horikawa, Eckart Viehweg, and to Professors Serge Cantat, Tien-Cuong Dinh for several valuable comments. {\it It is my great honor to dedicate this note to Professor Doctor Thomas Peternell on the occasion of his sixtieth birthday. He continues to inspire me through many interesting problems with his brilliant ideas since 1993}.

\end{document}